\documentclass[12pt]{article}
\usepackage{latexsym}
\usepackage{amsmath}

\begin{document}

\newcounter{thm}
\newtheorem{Def}[thm]{Definition}
\newtheorem{Thm}[thm]{Theorem}
\newtheorem{Lem}[thm]{Lemma}
\newtheorem{Cor}[thm]{Corollary}
\newtheorem{Prop}[thm]{Proposition}




                

\setlength{\baselineskip}{15pt}
\newcommand{\vlimsup}{\mathop{\overline{\lim}}}
\newcommand{\vliminf}{\mathop{\underline{\lim}}}
\newcommand{\Av}{\mathop{\mbox{Av}}}
\newcommand{\spec}{{\rm spec}}
\newcommand{\qed}{\hfill $\Box$ \\}
\def\textmc{\rm}
\def\({(\!(}
\def\){)\!)}
\def\R{{\bf R}}
\def\Z{{\bf Z}}
\def\N{{\bf N}}
\def\C{{\bf C}}
\def\T{{\bf T}}
\def\E{{\bf E}}
\def\H{{\bf H}}
\def\Prob{{\bf P}}
\def\M{{\cal M}}     
\def\F{{\cal F}}
\def\G{{\cal G}}
\def\D{{\cal D}}
\def\X{{\cal X}}
\def\A{{\cal A}}
\def\B{{\cal B}}
\def\L{{\cal L}}
\def\a{\alpha}
\def\b{\beta}
\def\e{\varepsilon}
\def\de{\delta}
\def\ga{\gamma}
\def\k{\kappa}
\def\la{\lambda}
\def\fa{\varphi}
\def\th{\theta}
\def\si{\sigma}
\def\t{\tau}
\def\om{\omega}
\def\De{\Delta}
\def\Ga{\Gamma}
\def\La{\Lambda}
\def\Om{\Omega}
\def\Th{\Theta}
\def\lan{\langle}
\def\ran{\rangle}
\def\lbr{\left(}
\def\rbr{\right)}
\def\const{\;\operatorname{const}} 
\def\dist{\operatorname{dist}} 
\def\Tr{\operatorname{Tr}}
\def\quadd{\qquad\qquad}
\def\n{\noindent}
\def\beq{\begin{eqnarray*}}
\def\eeq{\end{eqnarray*}}
\def\supp{\mbox{supp}}
\def\beqn{\begin{equation}}
\def\eeqn{\end{equation}}
\def\bp{{\bf p}}
\def\sg{{\rm sgn\,}}
\def\1{{\bf 1}}
\def\pf{{\it Proof.}}
\def\v2{\vskip2mm}
\def\n{\noindent}
\def\z{{\bf z}}
\def\x{{\bf x}}

\begin{center}
{\Large  Asymptotics of the densities    of \\the first passage time distributions
 for Bessel diffusions} \\
\vskip4mm
{K\^ohei UCHIYAMA} \\
\vskip2mm
{Department of Mathematics, Tokyo Institute of Technology} \\
{Oh-okayama, Meguro Tokyo 152-8551\\
e-mail: \,uchiyama@math.titech.ac.jp}
\end{center}

\vskip8mm

\begin{abstract}
 This paper concerns  the first passage times to a point $a >0$, denoted by $\sigma_a$,  of Bessel processes. We are interested in the case when  the process starts at $x>a$ and  compute the densities of the distributions of  $\sigma_a$  to obtain   the exact  asymptotic forms  of them  as  $t\to\infty$  that are valid uniformly in $x>a$ for every order of Bessel process.\footnote{
{\it key words}: first passage time,   exterior problem,  uniform estimate, Bessel diffusion\\
\n
}
 \end{abstract}
\vskip6mm

\section { Introduction and  main results}

This paper concerns   the first passage times to a point $a\geq 0$, denoted by $\sigma_a$,  of Bessel processes of order $\nu\in\R$. We are interested in the case when  the process starts at $x>a$ and compute the densities of the distributions of  $\sigma_a$  to obtain   the  exact asymptotic forms  of them  as   $t\to\infty$  that are valid uniformly in $x>a$ for each order $\nu$.
If $\nu=\pm 1/2$,  we have  well-known explicit expressions of them, which are  often used in various circumstances, while otherwise there has been  quite restricted information on them  until quite recently.   
 In the case when $0\leq x<a$ the distribution of $\sigma_a$ solves  a  boundary value problem of the associated  second order differential equation on the finite interval $(0,a)$ and the distribution of $\sigma_a$ or its density is  represented by means of  eigenfunction expansion (\cite{CT}, \cite{GS}, \cite{K}, etc.) and thereby we can obtain accurate estimates of them. In  the case $x>a$, however,   the region for the differential equation  is the infinite interval $(a, \infty)$  and   the corresponding representation is given by a Fourier-Bessel transform (cf. \cite{T1}: Section 4.10), which it seems not simple  a matter  to derive an asymptotic form of the density  directly from  and there have been only   a  few    partial  results as given in  \cite{S},    \cite{Ubh}, \cite{G} in which  $\nu=0$ or/and  relative ranges of $x$ are restricted at least for sharp estimates (in addition to  the cases $\nu =\pm 1/2$).  In the recent paper \cite{BMR}
  Byczkowski, Malecki and Ryznar  have computed an estimate  of the density for $\sigma_a$  for all values of $\nu$  by using  a certain integral representation  of it given in \cite{BR}: they obtain   upper and lower bounds  of the correct order of magnitude valid uniformly for all $t>0, x>a$, which  however does not give the exact  asymptotic form as we shall obtain in this paper (although in some cases their results are very close to and even finer than ours, see (i) of Remark 1 of the present paper).  
  Hamana and Matumoto  \cite{HM2}  have derived a similar (but,  in a significant point,  quite different) integral  representation  of  the  density of $\sigma_a$ for the case $x>a$ (as well as for the case $0\leq x<a$)  and computed
  an exact asymptotic form  of the density as $t\to\infty$  but with $x>a$ fixed.  
  
  The present investigation  is originally  motivated by a study of Wiener sausage of Brownian bridge in $\R^d$ joining the origin to a point $\x\in \R^d$ over a time interval  $[0,t]$
when $|\x|$ grows linearly with $t$ (\cite{Usaus2}).  The evaluation of the expected  volume of the sausage swept by a ball of radius $a$  can be reduced to that of the density for $\sigma_a$ with  arbitrary starting points $ >a$.  Although   only the case of the order $\nu =(d-2)/2$   is concerned there, the results for all orders $\nu \geq 0$  turns out to take parts in the evaluation.  
\v2

 Let $X^{\nu}_t$ be the Bessel process of order $\nu \in \R$, whose infinitesimal  generator $\L^{\nu}$ is given by
$$\L^{\nu} = \frac12 \frac{d^2}{d x^2} +\frac{2\nu+1}{2x}\cdot\frac{d}{d x} .$$
If $2\nu+2$ is a positive  integer, $X^{\nu}_t$ represents the radial part of the standard $(2\nu+2)$-dimensional Brownian motion.
If $\nu\geq -1$ we  write $d$ for $2\nu+2$:
 \beqn\label{d/nu}
 d=2\nu+2, ~~\mbox{or what is the same,}~~ 
 \nu =\frac{d}{2}-1;~~~
 \eeqn
 the process $X^{\nu}_t$ is sometimes called the $d$-dimensional Bessel process no matter whether $d$ is integer or not. 
Let $P_x$ be the probability law of the process $X^{\nu}_t$ started at $x\geq 0$ and $E_x$ the expectation by $P_x$.
Let $\sigma_a$ denote the first passage time of $X^{\nu}_t$ to $a>0$ and  $q^{\nu}(t,x;a) $ the density of the distribution of $\sigma_a$:
$$q^{\nu}(x,t;a) =\frac{d}{dt}P_x[\sigma_a \leq t].$$
We also write $q^{(d)}$ for $q^{\nu}$ where $d=2\nu +2$, if $\nu\geq -1$.

In what follows  we suppose $\nu\geq 0$ unless the contrary is stated  explicitly.  At the end of this introduction we shall observe that there is a simple   relation between  $q^{\nu}$ and $q^{-\nu}$ and  the case $\nu<0$  is reduced to the case $\nu>0$ and vice versa.   If  $\nu\geq 0$,  the origin is an entrance and non-exit boundary to the positive half line as is well-known. 
We shall use the two  indices $d$ and $\nu$ interchangeably,  understanding  that they are related by (\ref{d/nu}). 
Put 
$$ p_t^{\nu}(x)  = p_t^{(d)}(x)   = (2\pi t)^{-d/2} e^{-x^2/2t}.$$
For the process $X^\nu$  started at the origin $p_t^{(d)}(x) $ is the density of the distribution of $X^\nu_t$ w.r.t. the invariant measure $c_{d}x^{d-1}dx$, where $c_d= 2\pi^{d/2}/\Ga(\frac12 (d+1))$, the normalizing constant.  We prefer $q^{(d)}$ and $p_t^{(d)}$ to $q^\nu$ and $p^\nu_t$ and  in order to avoid confusion we shall  use the former notation  throughout  the paper except in a few occasions when the use of the latter one is definitely suitable.   

In \cite{Ubh} the present author obtains   the following result among others (see (\ref{q_2}) of Section 2 for another one). Put  $\k = 2e^{-2\ga}$, where  $\ga =-\int_0^{\infty} e^{-u}\lg u\, du$ (Euler's constant).
 
 \vskip2mm\n
 \begin{Thm}\label{thm1} (\cite{Ubh})~ If $\nu=0$, then  uniformly for $x>a$,  as $t\to\infty$
\begin{eqnarray}\label{q_x40}
q^{(2)}(x,t; a) =  \frac{\lg({\textstyle  \frac12} \k x^2/a^2)\,}{t(\lg (\k t/a^2))^2} e^{-\ x^2/2t} + \left\{\begin{array}{ll}   {\displaystyle \frac{2\ga \lg({t}/{x^2})}{t(\lg t)^3}+ O\bigg(\frac1{t(\lg t)^3} \bigg) }\quad & \mbox{for}~~~x^2<t,  \\[4mm]
{\displaystyle  O\bigg(\frac{1+[\lg (x^2/t)]^2}{x^2(\lg t)^3}\, \bigg) }\quad& \mbox{for}~~~x^2 \geq t.
\end{array}\right.
\end{eqnarray}
(Compare the first half with the corresponding one in Corollary 4 given below.)
\end{Thm}
\v2

This theorem does not identify  in any sense the asymptotic form   of $q^{(2)}(x,t;a)$  for $x>\sqrt{4 t\lg\lg t}$ (see Lemma \ref{lem3.1} in Section 3).  The objective of this  paper is to complement this ((\ref{R21}) in Theorem \ref{thm3} and Corollary  \ref{cor5}), and at the same time  also to obtain an asymptotic form  of $q^{\nu}$ for $\nu>0$ when the Bessel process is transient ((\ref{R2}) of Theorem \ref{thm3} and Theorem \ref{thm2}).

Define a function $\La_\nu(y), y\geq 0$ by 
$$\La_\nu(y)= \frac{(2\pi)^{\, \nu+1}}{2y^{\nu}K_{\nu}(y)}, \quad y>0$$
and $\La_\nu(0) = \lim_{y\downarrow 0}\La_\nu(y)$. Here $K_\nu$ is the modified Bessel function of second kind of order $\nu$; $\La_\nu(0)$ is well defined (see (\ref{K01}) below), and  
\[
\La_\nu(0)  = \frac{2\pi^{\nu+1}} {\Ga(\nu)} 
\quad\mbox{for}\quad\nu >0;\,\mbox{and}
\]
 $$\La_0(y) \sim \frac{\pi}{-\lg y}  \quad \mbox{as}\quad y \downarrow 0,$$ 
  in particular $\La_0(0)=0$. (By definition the ratio of two sides   of  $\sim $ tends to 1 in designated manner of taking limit.)   The main result of this paper is then  stated as follows.

\begin{Thm}\label{thm3}~  Uniformly for $x > a$, as $t\to\infty$, 
\beqn\label{R2}
q^{(d)}(x,t;a) =a^{2\nu}\La_\nu\bigg(\frac{a x}{t}\bigg) p^{(d)}_t(x) \bigg[1-\bigg(\frac{a}{x}\bigg)^{2\nu}\,\bigg]\Big(1+o(1)\Big)\qquad \mbox{if} \quad\nu>0
\eeqn
and  
\beqn\label{R21}
q^{(2)}(x,t;a) =p^{(2)}_t(x) 
\times
 \left\{ \begin{array} {ll}  
  {\displaystyle \frac{4\pi\lg(x/a)\,}{(\lg t)^2}\Big(1+o(1)\Big)    } \quad& (x  \leq \sqrt t\,),\\ [5mm]
 {\displaystyle  \La_0\bigg(\frac{a x}{t}\bigg) \Big(1+o(1)\Big)  }\quad&(  x   > \sqrt t\,).
  \end{array} \right.
\eeqn
If the  right-hand sides are multiplied by $e^{-a^2/2t}$, both the formulae (\ref{R2}) and (\ref{R21}) so modified  hold true  also as $x\to\infty$ uniformly for $t>0$.
\end{Thm} 

From the estimates of the density $q^{(d)}$ we can easily compute those of the distribution   $P_x[\sigma_a<t]$ and we shall carry out the computation that will be based on (\ref{R2}) and (\ref{R21}) in the last short section.
\v2\n
{\sc Remark 1.}~  (i)  In the recent paper \cite{BMR} 
  Byczkowski, Malecki and Ryznar  give estimates closely related to Theorem \ref{thm3}.  The main result of \cite{BMR}, their Theorem 2 (in its section 3),  is  a weaker version of the estimates (\ref{R2}) and (\ref{R21}):  the uniform  upper and lower  bounds of  the correct asymptotic order of magnitude  are obtained instead of  the exact asymptotic  form. For the case when  $x/t\to\infty$, however,  they derive  a very precise  estimate, finer than one given above, of which we present an explicit statement  shortly (Lemma \ref{thm5} below).   For $\nu>0$  their Proposition 5   identifies the asymptotic form of $q^{(d)}(x,t;a)$ in the case when $ x/t$ converges to a positive constant (as $t\to\infty$), the same  result as  included in Theorem  \ref{thm3} as a significant special case.   It is also noted that for each  $x>a$ fixed the formula (\ref{R2})   is given in \cite{HM2} but  with some coefficient being not explicit (see also \cite{BR}). The proofs in these papers 
rest on  certain integral representations of  $q^{(d)}(x,t)$ (given by \cite{BR} (in \cite{HM2}) and  by \cite{HM} (in \cite{BMR})) and the methods adopted therein  are quite different from ours.

(ii)  For every $\nu$, $\La_\nu(y)$ is an  increasing function of $y$.  Indeed, we have  
\beqn\label{K01}
 \La_\nu(y)= \frac{2\pi^{\nu+1}}{\int_0^\infty  \exp(-\frac1{4u}{y^2 })e^{-u } u^{\nu-1}du}~~~~~~~~~  (y>0),
\eeqn
as is readily deduced from the identity  $K_\nu(z) =\frac12 (z/2)^{\nu}\int_0^\infty e^{-\frac1{4u}{z^2 } -u} u^{-\nu-1}du
$ $(|\arg z|<\frac{\pi}{4})$  (\cite{L}, p.119),  of which   the  right-hand side is invariant under the replacement of $\nu$ by $-\nu$. 
It is also noted that for each $\nu\geq 0$, $K_\nu(y) =\sqrt{\pi/2y} \, e^{-y} ( 1+O(1/y))\, (y\to\infty)$, so that
\beqn\label{Kinfty}
\La_\nu(y) =(2\pi)^{\nu+1/2} y^{-\nu+1/2}\, e^{y} ( 1+O(1/y))~~~~(y\to\infty).
\eeqn

\v2

Detailed   estimation of $q^{(d)}(x,t;a) $ inside the parabolic regions  $x^2<C t$ $ (C>1)$ would be of fundamental importance. We next give an extension of Theorem 1 to $\nu>0$ in this respect, which partially  prepares for  the proof of Theorem \ref{thm3} in an 
  obvious way (see (iii) of Remark 2 below).  
  
  For $\nu>0$  we have the Green function, $G^{(d)}(x,y)$ say; we need to bring in $G(x)=G^{(d)}(x) ;= G^{(d)}(x,0)= G^{(d)}(0,x)$, or explicitly 
  $$
G(x) =\int_0^\infty p^{(d)}_t(x)dt = \frac{\Ga(\frac{d}2 -1) }{2 \pi^{d/2}}\cdot \frac1{x^{d-2}}.
$$ 
The next theorem (a reduced version of Propositions \ref{P1} and \ref{P2}  in Section 2) gives a fairly fine estimate of $q^{(d)}$ for $\nu>0$ in the case when  $x\leq \sqrt{2(\nu\wedge 1)t\lg t}$ (as in Theorem \ref{thm1} in a sense). 
 Here as well as in what follows 
  $c\vee b= \max \{c,b\},  c\wedge b=\min\{c,b\}$ for $c, b$ real.
 It is noted that $a^{2\nu}\La_\nu(0)=1/G(a)$. 

\vskip2mm 
   
 \begin{Thm} \label{thm2}~ Let $\nu > 0$. If $\nu\neq 1$, then uniformly for $x>a$, as $t\to\infty$
 \begin{eqnarray}\label{R1a}
 q^{(d)}(x,t;a) &=& \frac{(a^2/2)^\nu}{\Ga(\nu)\, t^{\nu+1}}  \Bigg[ e^{-x^2/2t}- \bigg(\frac{a}{x}\bigg)^{2\nu} e^{-a^2/2t}\,\Bigg] +O\Bigg(\frac{1\wedge ({\sqrt t}/{x})^{\nu+(\frac12 \wedge \nu)}}{t^{\nu+1+(\nu\wedge 1)}}\Bigg)\\
  &=&\frac{1}{G(a)} \,\Bigg[ p_t^{(d)}(x)- \bigg(\frac{a}{x}\bigg)^{d-2} p_t^{(d)}(a)\,\Bigg] +O\Bigg(\frac{1\wedge ({\sqrt t}/{x})^{[ (d-1)/2] \wedge (d-2)}}{t^{d/2 +(\nu\wedge 1)}}\Bigg).\label{R1b}
  \end{eqnarray} 
In the case  $\nu =1$ the same estimate holds true if the error term given by the $O$-symbol is replaced by
$$O\bigg(\frac{1+\lg (t/x^2)}{t^3}\bigg)~~~\mbox{ for}~~
 a<x<\sqrt t;~~~~ O\bigg(\frac{\lg t}{t^3} \bigg(\frac{\sqrt t\,}{x}\bigg)^{3/2}\bigg)~~~\mbox{for}~~
 x>\sqrt t.$$
  \end{Thm}
\vskip2mm\n
 {\sc Remark} 2.~ (i) ~For  random walks on the $d$-dimensional square lattice $\Z^d$ we have analogues of Theorems \ref{thm1} and \ref{thm2} \cite{U2}.  The form of the principal  term in the formula (\ref{R1b}) is  intrinsically the same as and in fact suggested by that corresponding to  the walks.  For $d=2$ Theorem \ref{thm2} (or Proposition \ref{lem3.0})  provides 
 an  improvement  of Theorem 1.4 of \cite{U2} in view of Theorem 1.5 of it. 
 
  (ii) ~In the proof of Theorem \ref{thm2} we give a more precise expression of the error term, which shows that its order of magnitude cannot be improved at least for $x<\sqrt t$.

(iii)~ If the range is restricted to $2a<x<\sqrt t$,  (\ref{R2}) is  immediate consequence  of Theorem \ref{thm2};    
the case $a<x\leq  2a$ is contained in Proposition \ref{P1}, which also implies  (\ref{R2}).  In view of these results for $x<\sqrt t$  the essential ingredient of (\ref{R2}) is now the estimate in the region  $x>\sqrt t$.  With certain  additional results employed the same can be said  for $\nu=0$ (for details  see {\sc Remark}  given  at the end of Section 2).

\v2
Two dimensional case is particularly interesting and deserves to  be described here  in more detail. Restricting to the region $x>\sqrt t$  we may state the formula (\ref{R21})  as follows: 
  uniformly for $ x > \sqrt t$,  as $t\to\infty$ 
 \beqn\label{EQ4}
 q^{(2)}(x,t, a)  =   \frac{1}{2K_0(a x/t)}\cdot\frac{e^{-x^2/2t}}{t}\Big(1 + o(1)\Big).  
 \eeqn
Substitution from the formulae (\ref{Kinfty})  and $K_0(u) =- \lg (\frac12 e^\ga u) + O(u^2\lg u)$ ($u\downarrow 0$)  makes the right-hand side above explicit if $x/t$ goes to $0$ or $\infty$.  We shall actually compute   errors  in  the formula  (\ref{EQ4}) (see Proposition \ref{lem3.0} and Lemma \ref{thm5}).  Taking account of   these comments  the next result is essentially a corollary of  Theorem \ref{thm1} and the proof of  Theorem \ref{thm3}.
\begin{Cor}\label{cor5} ~Let $\nu=0$. Uniformly in $x>a$, as $t\to\infty$, 
\begin{eqnarray}\label{q1} 
{q^{(2)}(x,t, a)}  &=&   \frac{2\lg (x/a)}{ (\lg  (t/a^2))^2}\cdot\frac{e^{-x^2/2t}}{t}\bigg[1 +O\bigg( \frac{1}{ \lg t}\bigg) \bigg]
 ~~~~~~~~~~~~~\mbox{if}~~~~ ~a<x \leq \sqrt t, \nonumber \\[1mm]
&=& \frac{1}{ 2\lg (t/ ax)} \cdot\frac{e^{-x^2/2t}}{t}\bigg[1 +O\bigg( \frac{1}{ \lg (t/x)}\bigg) \bigg]
~~~~~~~~~\mbox{if}~~~~ x/t \to 0, ~ x >\sqrt t, \nonumber \\[1mm]
&=& \sqrt{\frac{ ax}{2\pi  t}} \cdot\frac{e^{-(x-a)^2/2t}}{t}\bigg[1 +O\bigg( \frac{t}{ x}\bigg) \bigg]
 ~~~~~~~~~~~~~~ ~\mbox{if}~~~~x/t \to \infty.
\end{eqnarray}
\end{Cor}

The factor $e^{-(x-a)^2/2t}$ in the last formula (\ref{q1}),  asymptotically  equivalent to $e^{-x^2/2t}$ when $x/t\to 0$,  may be understood  to  be natural by comparing with  the Gaussian kernel $p^{(2)}_t(x-a)$ (see also (\ref{33}) in Section 2).

The next result from    \cite{BMR} (Lemma 4) gives a fine estimate  
in the case when  $x/t\to\infty$.  It in particular gives asymptotics of  $q^{(d)}$ as $t\downarrow 0$; also shows  that the dependence on $d$ of the leading term   in this case comes only from  the factor $(a/x)^{(d-1)/2}$. 
\begin{Lem}\label{thm5} (\cite{BMR})~  For each $\nu\geq 0$ it holds that  uniformly
 for  all $t>0$ and $x>a$, 
\beqn\label{result6}
 q^{(d)}(x,t;a) =  \frac{x-a}{\sqrt{2\pi}\, t^{3/2}} e^{-(x-a)^2/2t}\bigg(\frac{a}{x}\bigg)^{(d-1)/2}\Bigg[1+ \frac{\b t}{ ax}\bigg(1+O\bigg(\sqrt t \wedge  \frac{t}{x-a}\bigg) \bigg)\Bigg].
 \eeqn
 where $\b = (d-1)(3-d)/8 = (\frac14- \nu^2)/2$.
 \end{Lem}
 
\v2
From the scaling property of Bessel processes  it follows that
$$q(x,t;a) = a^{-2} q(x/a, t/a^2;1).$$
 For the proofs of  the foregoing theorems we shall mostly  consider only the case $a=1$  and write $q(x,t)$ for $q(x,t;1)$.

\v2
The estimation of $q(x,t)$ will be  made  in the following  three cases
\v2
 (i) ~~ $x< \sqrt t$; ~~(ii) ~~$\sqrt t < x \leq Mt$ (with $M$ arbitrarily fixed) ; ~~(iii) ~~~ $x/t\to \infty$,
 \v2\n
of which  the cases (i) and (ii)   will be  discussed in Sections 2 and 3, respectively.  The methods employed in these cases are different from one another. Roughly speaking, for  the case (i) the estimation is based on the well known formula for the Laplace transform of $q^{(d)}(x,\cdot)$, to which we apply the Laplace inversion formula; some computation using  the Cauchy integral theorem then  leads to  somewhat finer estimates (Proposition  \ref{P1}) than those given in  Theorem \ref{thm2}. 
For the case   (ii) we exploit   the fact that any  Bessel process of order $\nu> -1$  can be decomposed as a sum  of two  independent Bessel processes and apply the result of the case (i). This gives some error estimate to the asserted  asymptotic form  of $q$ in the case $x/t\to 0$ (Propositions \ref{lem3.0} and \ref{thm3.2}).  To include the case $x/t \to v>0$ an additional argument is employed.

For (iii)  Lemma \ref{thm5} provides a better estimate than required for Theorem \ref{thm3}. The proof of Lemma \ref{thm5} rests on the integral representation obtained in \cite{BR} and the derivation from it is involved. A relatively easier  proof for the relevant estimate in  Theorem \ref{thm3}  can be provided by making use of
 the following  probabilistic expression: 
\beqn\label{Eq0}
q^{(d)}(x,t;a) = \frac{x-a}{\sqrt{2\pi t^3}} e^{- (x-a)^2/2t}\bigg(\frac{a}{x}\bigg)^{(d-1)/2} E^{BM}_x\Bigg[\exp\bigg\{ \b\int_0^{t} \frac{ds}{B_s^{2}}\bigg\}\,\Bigg| \sigma_a =t \Bigg],
\eeqn
where  the conditional expectation is taken w.r.t. the probability measure of the standard linear Brownian motion $B_t$.  (This identity is readily derived from the well known formula for $q^{(1)}$ by   using the Cameron-Martin-Girsanov formula.)     In the case $x/t\to\infty$, $t>1$  the   estimate    (\ref{result6})  but  with the second term   in the big square brackets in it being  replaced by a  less exact $O(t/x)$  can be derived   from  (\ref{Eq0}) by some comparison argument based on   the diffusion  equation that is associated with the expectation in (\ref{Eq0}) via  the  Kac formula (cf.  \cite{Usaus2}), of which  we will not go into further details. By the way it is noted that
the expression of  $q^{(d)}$ above  verifies, on expanding the exponential, the estimate (\ref{result6}) of the case  $t\downarrow 0$   with the same replacement for  the second term  as above.

 If $x/t$ is large enough, one can evaluate the conditional expectation in  (\ref{Eq0}) for large $t$ directly as mentioned right above, of which  the dependence on $\nu$ comes only from $\b$. Otherwise, however,  a direct evaluation of it seems hard.  Our results on $q^{(d)}$ rather give 
a precise estimate of it valid uniformly in  $x$, which turns out to be  useful: in \cite{Usaus2} we exploit the estimate to derive an asymptotic form of the  space-time distribution  of the hitting  of a ball by $d$-dimensional Brownian motion.

Throughout the paper $C, C', C'',$ etc.  will be used to denote constants whose precise values are not important 
for the present purpose; the same letter  may designate   different constants  depending on the occasions where it occurs.

 \v2
 We conclude the present section  by mentioning  some simple facts for the case $\nu<0$.  The probability that the Bessel process  $X_t^{|\nu|}$ of order $|\nu|$ started at $x>a$ hits  $a$ in a finite time is given by $h(x)=(a/x)^{2|\nu|}$, which is a harmonic function for the process restricted on $[a, \infty)$ with killing at $a$, and   the conditional process conditioned on this event 
 is a $h$-transform of it. On identifying the generator   this conditional process is a  Bessel process of order $\nu\,(<0)$. Hence  
$$q^\nu(x,t;a) =  q^{|\nu|} (x,t;a)\frac{x^{2|\nu|}}{a^{2|\nu|}}.$$
(This also  follows from (\ref{T(2)})  below.)
 Every Bessel process of negative order visits  the origin in a finite time with probability one and we have explicit formula
 \beqn\label{EQ11}
 q^\nu(x,t;0) = \La_{|\nu|}(0)x^{2|\nu|}p_t^{|\nu|}(x)
 \eeqn
 (for the derivation let $a\downarrow 0$ in (\ref{T(2)}) and use (\ref{2.11}), both given in Section 2).
 By a comparison argument we have the inequality
 $$\int_0^t q^\nu(x-a,s;0)ds  <\int_0^t  q^\nu(x,s;a)ds ~~~(x>a, t>0) ~~~~~\mbox{if} ~~~ -2^{-1} <\nu<0$$
 and the same  one but  in the  opposite direction  if  $\nu<-2^{-1}$. 
 

  \section{Proof of Theorem \ref{thm2}}
  
For any  $\nu\in \R$, 
\beqn\label{T(2)}
E_x[\exp\{-\la \sigma_a\}] 
=\frac{K_{\nu}(x\sqrt{2\la})x^{-\nu}}{K_{\nu}(a\sqrt{2\la})a^{-\nu}}~~~~~(\la>0, x>a>0),
\eeqn
as is well-known and may be derived  by solving the problem: $\L^{\nu}U=\la U$ $(x > a)$ with the lateral conditions  $U(a+0)=1$ and $U$ being  positive and decreasing, of which solution  is unique. (See (\ref{LT}) below.)

For $\nu=1/2$  (i.e., $d=3$) we have a particularly simple expression of $q^{(3)}$:  for $x>a$
\beqn\label{33}
q^{(d=3)}(x,t;a)=\frac{ae^{-(x-a)^2/2t}}{t\sqrt {2\pi t}}\bigg(1-\frac{a}{x}\bigg),
\eeqn
which trivializes this special  case of  Theorems \ref{thm2} and \ref{thm3} and is helpful for making a guess at  the asymptotic form of $q^{(d)}$  in  general cases.  

In what follows  we  let $\nu\geq 0$  and, when  there is no risk of confusion,  we suppress the super-script$\,^{(d)}$ from  $q^{(d)}(x,t;a) $ and  $p_t^{(d)}(x) $ except for the statement of propositions or lemmas.  

Put
 $G_\la(x) =\int_0^\infty p_t(x) e^{-\la t}dt$. We know 
\beqn\label{2.11}
G_\la(x)  =\frac2{(2\pi)^{d/2}}\bigg(\frac{\sqrt {2\la}}{x}\bigg)^\nu K_\nu(x\sqrt{2\la})
\eeqn
 (\cite{E}, p.146). 
It is convenient (and natural) to  write the representation (\ref{T(2)}) in the form 
\beqn\label{LT}
E_x[\exp\{-\la \sigma_a\}] =\frac{G_\la(x)}{G_\la(a) },
\eeqn
For $\nu>0$ let $G(x)  =\lim_{\la\downarrow 0}G_{\la}(x)$, so that
$$
G(x) =\int_0^\infty p_t(x)dt = \frac {\Ga(\frac{d}2 -1) }{2 \pi^{d/2} x^{d-2}}.
$$
If $\nu$ is not an  integer, $K_\nu(z) = \Big(\pi /2\sin (\pi \nu)\Big)\Big[I_{-\nu}(z)- I_\nu(z)\Big]$, where $I_\nu$ is the modified Bessel function of the first kind of order $\nu$  and given by 
\beqn\label{E_of_I}
I_\nu(z)= \bigg(\frac{z}{2}\bigg)^\nu \sum_{k=0}^\infty \frac{(z/2)^{2k}}{\Ga(\nu+k+1)\Ga(k+1)}
\eeqn
for $|\arg z|<\pi$ (\cite{L}, p.108).

\begin{Prop}\label{P1} ~  Let $\nu>0$ and $M>0$. For $a<x< M\sqrt t,$
\beqn\label{EqP1}
q^{(d)}(x,t;a) = \frac1{G(a)}\bigg[p^{(d)}_t (x) -\bigg(\frac{a}{x}\bigg)^{d-2}p^{(d)}_t(a)\bigg] + a^{-2} \eta(x/a,t/a^2),
\eeqn
with the function  $\eta(x,t)$, $x>1, t>2$ admitting the estimate
\beq
&&\eta(x,t) = O\bigg(\frac{1-x^{-1}} {t^{\nu+2}}\bigg) ~~~~\mbox{if} ~~\nu>1; ~~= O\bigg(\frac{1-x^{-1}} {t^{2\nu+1}}\bigg) ~~~~\mbox{if} ~~0<\nu <1 \\
&&\mbox{and} ~~~~
\eta(x,t) = O\bigg(\frac{1-x^{-1} } {t^{\nu+2}}\, \lg t\bigg) ~~~~\mbox{if} ~~\nu=1.
\eeq
(See (\ref{P_eq5}), (\ref{P_eq61}) and (\ref{P_eq62}) for more exact forms of $\eta(x,t)$.) 
\end{Prop} 
\v2\n
{\sc Remark} 3. ~One might suspect  that the function  $q^*(x,t):= [G(a)]^{-1}p^{(d)}_t(x-a)[1- (a/x)^{d-2}]$, an analogue to the explicit  form  of $q^{(3)}$, can take  the  place of  the leading term in the formula  (\ref{EqP1}). Since
the difference of them  is at most  the magnitude of  $O(x^{(2-2\nu)\vee 1}/t^{\nu + 2})$, this is true  if $\nu < 1$; in the case $\nu>1$, however,   the difference becomes much larger than $\eta(t,x)$ as $x$ gets large, so
that the replacement causes  a larger error term.
\v2\n
{\it Proof of Proposition \ref{P1}.} ~ 
The Laplace inversion of (\ref{LT})  gives
\beqn \label{P1.1}
q(x,t;a) =  \frac1{2\pi i}\int_{-i\infty}^{i\infty} \frac{G_{z}(x)}{G_{z}(a)} e^{tz}dz.
\eeqn
Here and  in what follows   the functions $z^{\nu/2}$ and $\lg z$ involved in the integrand (see (\ref{E_of_I}), (\ref{A})) take their principal values. 
For the evaluation of the integral above we follow the argument made in \cite{U2} for the random walk of dimensions $d\geq 3$.  Motivated by it we decompose
\begin{eqnarray}
\frac{G_{\la}(x)}{G_{\la}(a)}&=&\frac{G_{\la}(x) - (a/x)^{2\nu} G_{\la}(a) }{G_{\la}(a)}  + \frac{a^{2\nu}}{x^{2\nu}} \nonumber\\
&=& \frac{G_{\la}(x) - (a/x)^{2\nu} G_{\la}(a) }{G(a)}  + R(\la; x) 
\label{P1.3}
\end{eqnarray}
where
$$R(\la; x)= \bigg[  \frac{1}{G_{\la}(a)}   -\frac{1}{G(a)}\bigg] \Big[G_{\la}(x) - (a/x)^{2\nu} G_{\la}(a) \Big] + \frac{a^{2\nu}}{x^{2\nu}}. $$
 By the definition of $G_\la(x)$  the contribution to (\ref{P1.1}) of the first term on the right-hand side of (\ref{P1.3}) equals 
$$
\frac1{G(a)}\Big[ p_t(x) -(a/x)^{2\nu} p_t(a)\Big].
$$
The error term in Proposition \ref{P1}  is then written as
\beqn\label{eta}
\eta(x,t) =\frac1{2\pi i}\int_{-i\infty}^{i\infty} R(z,x) e^{tz}dz.
\eeqn

 As $z\to \infty$
\beqn\label{K1}
K_\nu(z) = (\pi/2 z)^{1/2}e^{-z}(1+O(1/z))~~~~~ (|\arg z|\leq \pi).
\eeqn 
Hence $G_{z}(x)/G_{z}(a) = O(e^{-(x-a)\sqrt{2z}})$   as $z\to \infty$ (with $x$ fixed) and we see that the function $R(z,x)$ rapidly approaches  zero as $z\to \infty$ in the sector $|\arg z|>\pi -\de$ for any $\de>0$. This together with the fact  that $K_\nu(z)$ has no zeros on the right half plane $\Re z \geq 0$ (cf. \cite{W}: p. 511) (hence $G_z(x)$ has no zeros on $-\pi\leq \arg z\leq \pi$)  permits us to apply  Cauchy's integral formula to transform (\ref{eta}) into
\beqn \label{P_eq0}
\eta(x,t)
 =  \frac{1}{2\pi i} \int_{0}^\infty \Big [-R(-u+i0,  x)+ R(-u-i0,x)\Big]e^{-tu}du,
\eeqn
where $R(-u+i0,  x) = \lim_{y\downarrow 0} R(-u+iy,  x)$ and analogously for $R(-u-i0,x)$.
It is noted that  the integral in (\ref{P_eq0}) is not affected by subtraction  of any  entire function from $R$.
\v2

{\it Let $\nu$ be not an integer. }~ We show that uniformly for  $a<x<M\sqrt t$,
\begin{eqnarray}\label{P_eq5}
\eta(x,t) =  \frac{ -\nu-1}{2^{\nu+1}(\nu-1)\Ga(\nu)} \Bigg(1-\frac{1}{x^{2\nu}} +\frac 
{x^2-1}{x^{2\nu}}\Bigg) \frac1{t^{\nu+2}}  
 - B^2_0\frac{1-x^{-2\nu}}{\Ga(-2\nu)t^{2\nu+1}} + r(t,x)
\end{eqnarray}
with
 \beq
|r(t,x)|& \leq&  \frac{C}{t^{2\nu+1}} \Bigg[\bigg(\frac{x^2}{t}\bigg)\vee \frac{1}{t^\nu} \Bigg] \bigg(1-\frac1{x}\bigg)~~~~~(0<\nu<1),
\\
&\leq &   \frac{C}{t^{\nu+2}} \cdot\frac{x^2}{t} \bigg(1-\frac1{x}\bigg)~~~~~~~~~~~~~~~~~~~~(\nu >1),
\eeq
where the constant $B_0$ is given  in  (\ref{AB})   below.
 The estimation of $\eta$  is simple apart from the uniformity in $x$,  which we must take care of  in dealing with the dependence on $x$. 
Let  $a=1$ for simplicity. Recall  the definitions of  $G_\la(1)$ and $G(1)$.   From   the power series expansion  of $I_\nu$ given in (\ref{E_of_I})   we then deduce
straightforwardly  
\begin{eqnarray}\label{P_eq6}
\frac{G_\la(1)}{G(1)}&=& \frac{(\sqrt{2\la})^{\nu}[ \,I_{-\nu}(\sqrt{2\la}\,)- I_\nu(\sqrt{2\la}\,)]}{\lim_{\la\downarrow 0}(\sqrt{2\la})^{\nu}I_{-\nu}(\sqrt{2\la}\,)}  \nonumber\\
&=&\Ga(1-\nu)\Bigg(\sum_{k=0}^\infty \frac{(\la/2)^k}{\Ga(-\nu+k+1)k! } + \bigg(\frac{\sqrt{2\la}}{2}\bigg)^{2\nu}\sum_{k=0}^\infty \frac{(\la/2)^k}{\Ga(\nu+k+1) k! } \Bigg) \nonumber\\
& =& \Big[1 + A_1\la +A_2 \la^2+\cdots\Big] - \la^{\nu}\Big[B_0+B_1\la+\cdots\Big]
\end{eqnarray}
with
\beqn\label{AB}
A_1 = [2(1-\nu)]^{-1}, ~~~~B_0= 2^{-\nu}\Ga(1-\nu)/\Ga(\nu+1).
\eeqn
Let  $f(\la)=  G_\la(1)/G(1) - 1$. Then
\beqn\label{P_eq1} \frac{1}{G_{\la}(1)}   -\frac{1}{G(1)} = \frac{1}{G(1)}\bigg[\frac1{1+f(\la)} -1\bigg] = - \frac{ f(\la) }{G(1)}+ \frac{[f(\la)]^{2}}{ G_{\la}(1)}.
\eeqn
Also, noting  
$$x^{2\nu}G_\la(x)=G_{x^2\la}(1),$$
we obtain
\begin{eqnarray}\label{P_eq2}
x^{2\nu}\times \frac{G_{\la}(x) - (1/x)^{2\nu} G_{\la}(1) }{   G(1)} &=& \frac{ G_{x^2\la} (1)-  G_{\la}(1) }{G(1)}\nonumber\\
 &=& A_1\la  (x^2-1) - B_0  \la^\nu  (x^{2\nu}- 1)+ H_\la(x),
\end{eqnarray}
where $H_\la(x)= \Big[A_2 \la^2(x^4 -1)+\cdots\Big] +  \Big[ B_1\la^{\nu+1} (x^{2\nu+2}- 1)+\cdots\Big] $, the remainder term.   Since $f(z)=O(|z|^{\nu/2-1/4})$ as  $z\to \infty$ along the negative real line, the equality   (\ref{P_eq6}) entails 
  $$ f(\la)= A_1\la  - B_0 \la^\nu + O(|\la|^{(\nu\wedge 1)+1})~~~(|\la|<1)~~\mbox{and}~~
|f(-u\pm i0)|\leq C|u|^{\nu/2}~~(u>0).$$ 
From the power series expansions of  $I_{\pm \nu}(z)$ it follows that $A_k = O(1/\Ga(-\nu+k+1)k!),~ B_k =O(1/\Ga(\nu+k+1) k! ).$

With these preliminary discussions  we now  compute the integral in (\ref{P_eq0}). Using (\ref{P_eq1}) and (\ref{P_eq2}) we make decomposition
\beqn\label{EQ28}
x^{2\nu}R(\la;x) = - f(\la)\Big ( A_1\la (x^2-1) -B_0 \la^\nu  (x^{2\nu}- 1)\Big)+T_1(\la)+ f(\la) H_\la(x) +1,
\eeqn
where $T_1(\la) =[f(\la)]^2[G_{x^2\la}(1)- G_\la(1)]/G_\la(1); $  the principal  part would be involved in the first term  in view of (\ref{P1.3}), (\ref{AB}) and (\ref{P_eq2}).

First  consider the contribution of $T_1$  and observe that  for  $u>0$,
\beq
|T_1(-u+i0)-T_1(-u-i0)| &\leq& C_0u^{2(\nu\wedge 1)}\Big[u^\nu(x^{2\nu}-1) +C_1u^{\nu+1}(x^{2\nu+2}-1)+\cdots\Big] \\
&&+~ C_0'u^{2\nu}\Big[ u(x^2-1) + C_1' u^2(x^4-1)+\cdots\Big]
\eeq
with certain constants $C_k, C_k'$ that are dominated by a constant multiple of $2^k k^{3\nu}/(k!)^2$. Here we have exploited the fact that the terms of integral powers $c_n\la^n$ involved in $T_1(\la)$ cancels out in the difference on the left-hand side.
Employing the simple inequality  $x^s -1\leq (1\vee s)(1-x^{-1})x^s$ valid  for all $x>1, s>0$ we infer that for $1<x<M\sqrt t$,
\begin{eqnarray}\label{P_eq3}
&&\int_0^\infty \Big|T_1(-u + i0)- T_1(-u - i0)\Big|e^{-tu}du \nonumber \\
&&\leq  C\bigg(1-\frac1{x}\bigg) \Bigg(\frac{1}{t^{2(\nu\wedge 1)+1}} \sum_{k=0}^\infty \frac{k^{3\nu+1}4^k x^{2(\nu+k)}}{t^{\nu+k} k!} + 
\frac1{t^{2\nu+1}}\sum_{k=1}^\infty \frac{k^{3\nu+1} 4^k x^{2k}}{t^k k!}   \Bigg)\nonumber\\
&&\leq C' \bigg(1-\frac1{x}\bigg)\Bigg[\frac{1}{t^{2(\nu\wedge 1)+1}} \bigg(\frac{x^{2}}{t}\bigg)^\nu + 
\frac1{t^{2\nu+1}}\bigg(\frac{x^{2}}{t}\bigg)   \Bigg],
\end{eqnarray}
where $C'$  depends on $M$.

Secondly, in the same way we see that if $T_2(\la)= f(\la) H_\la(x)$, 
\begin{eqnarray}\label{P_eq4}
\int_0^\infty  |T_2(-u + i0) -T_2(-u  - i0)|e^{-tu}du 
\leq C \bigg(1-\frac1{x}\bigg)\bigg(\frac{x^{2}}{t}\bigg)^{\nu+1}  
\frac1{t^{\nu\wedge 1+1}}.
\end{eqnarray}

Thirdly  putting
$$F(\la)= - f(\la)\Big ( A_1\la (x^2-1) -B_0 \la^\nu  (x^{2\nu}- 1)\Big),$$
we have
$$
  F(\la)= A_1 B_0 \la^{\nu+1}   \Big(x^{2\nu}-1+{x^2-1}\Big)(1+ C_1\la +\cdots) -  B_0^2\la^{2\nu}(x^{2\nu}-1)(1 +C'_1\la+\cdots)
$$
apart from the difference by an   entire function,  and for $s>-1$,
$$\frac1{2\pi i}\int_0^\infty\Big[ -(-u+i0)^s +(-u-i0)^s\Big]e^{-tu}du = \frac{1~}{\Ga(-s)t^{s+1}}.$$
(Here the identity  $\Ga(1+s)\sin \pi s =-\pi / \Ga(-s)$ is used;  remember that $1/\Ga(-n) =0$ if $n$ is non-negative integer.) Hence 
\begin{eqnarray}
&&\frac1{2\pi i}\int_0^\infty  [-F(-u+i0) +F(-u-i0)] e^{-tu}du  \nonumber \\
&&= \Bigg( A_1 B_0\frac{x^{2\nu}-1+{x^2-1}}{\Ga(-\nu-1)t^{\nu+2}}  -B^2_0\frac{x^{2\nu}-1}{\Ga(-2\nu)t^{2\nu+1}} \Bigg)
 \Big[1+ O(1/t)\Big]
 \label{F_exp}.
\end{eqnarray}

Finally collecting the bounds (\ref{P_eq0}) and (\ref{P_eq6}) through (\ref{P_eq4}) (of which we divide each formula by $x^{2\nu}$ since we have multiplied it in (\ref{P_eq2}) and (\ref{EQ28})), noting $B_0/\Ga(-\nu-1) =(\nu+1)/2^\nu\Ga(\nu)$  and making elementary comparison of terms that appear on the right-hand sides of them  we find (\ref{P_eq5}) to be true. 
\v2

{\it Let $\nu$ be a positive integer.}~ The arguments are similar to the above. In place of (\ref{P_eq6}) we have
\begin{eqnarray}\label{A}
\frac{G_\la(1)}{G(1)}&=&\frac{2}{(\nu-1)!}\left[\,\sum_{k=0}^{\nu-1} \frac{(\nu-k-1)!}{ 2(k!)} \bigg(\frac{-\la}{2}\bigg)^k
-\bigg( \frac{-\la}{2}\bigg)^{\nu}\Bigg(\ga+\lg \sqrt{\frac{\la}{2}}\, \Bigg)\sum_{k=0}^\infty \frac{(\la/2)^k}{k! (\nu+k)!} \,\right]  \nonumber\\
&& + ~\bigg( \frac{-\la}{2}\bigg)^{\nu} g(\la) \nonumber\\
&=& \Big[1 + A_1\la +A_2 \la^2+\cdots\Big] - \la^{\nu}(\lg \la)\Big[B_0+B_1\la+\cdots\Big], 
\end{eqnarray}
where
  $g(\la) =a_0+a_1\la +\cdots$ is a certain  entire function (with $a_0=1$ for $\nu=1$), 
and 
\begin{eqnarray}
\label{AB2}
A_1 = \left\{\begin{array} {ll} -2^{-1}[1+\lg(2e^{-2\ga}) ] ~~& (\nu=1),\\
-[2(\nu-1)]^{-1}  ~&(\nu\geq 2) \end{array}\right. ~~\mbox{ and}~~~  B_0= \frac{(-1)^\nu}{2^\nu \nu!(\nu-1)! }.  
\end{eqnarray}
Noting 
 $$(2\pi)^{-1}\Im\Big(-[\lg (-u+i0)]^k+[\lg (-u-i0)]^k \Big) = -k(\lg u)^{k-1}~~~~ ~(k=1, 2, u>0)$$ we compute the integral  in (\ref{P_eq0}) to
  see that  if $\nu$ is an integer greater than 1, then
\begin{eqnarray}\label{P_eq61}
\eta(x,t) = - (-1)^{\nu+1}(\nu+1)!   A_1 B_0 \Bigg(1-\frac{1}{x^{2\nu}} +\frac 
{x^2-1}{x^{2\nu}}\Bigg) \frac1{t^{\nu+2}}  + O\Bigg(\frac{1-x^{-1}}{ t^{\nu+2}}\cdot \frac{x^2}{t}
\Bigg).
\end{eqnarray}
(The coefficient of the leading term  coincides with one appearing in   (\ref{P_eq5}) [since $1/\Ga(-n) =0$ if $n$ is a non-negative integer], so that  (\ref{P_eq5}) is valid for $\nu =2, 3, 4,\ldots$.)
If $\nu=1$ and $f$ is defined as before, the leading term of $- f(\la)\Big[ G_{x^2\la}(1)-G_\la(1)\Big]/G(1)$ being $$(B_0^2\la^2 \lg  \la)(- x^2\lg (x^2\la) +\lg \la),$$
 we see 
\begin{eqnarray}\label{P_eq62}
\eta(x,t) = - 4 B_0^2 \frac{ (1-x^{-2})\lg t - \lg x} {t^{\nu+2}} 
+ ~O\Bigg(\frac{1-x^{-1}}{t^{\nu+2}}\Bigg).
\end{eqnarray}
The foregoing two formulae obviously imply the desired bounds for $\eta(x,t)$.
The proof of Proposition \ref{P1} is complete.   \qed

We still need to obtain error estimates for  $|x|>\sqrt t$.

\begin{Prop}\label{P2} ~  Let $\nu>0$. Then the function  $\eta(t,x)$  defined via (\ref{EqP1}) 
enjoys  the estimates 
\beq
&&\eta(t,x) = O\Bigg(\frac{1} {t^{\nu+2}}\bigg(\frac{\sqrt t\,}{x}\bigg)^{\nu+\frac12}\Bigg) ~~~~\mbox{if} ~~\nu>1; ~~= O\Bigg(\frac{1} {t^{2\nu+1}} \bigg(\frac{\sqrt t\,}{x}\bigg)^{\nu+(\frac12 \wedge \nu)}\Bigg) ~~~~\mbox{if} ~~0<\nu <1 \\
&&\mbox{and} ~~~~
\eta(t,x) = O\Bigg(\frac{\lg t} {t^{\nu+2}}\, \bigg(\frac{\sqrt t\,}{x}\bigg)^{\nu+\frac12 }\Bigg) ~~~~\mbox{if} ~~\nu=1.
\eeq
that are valid uniformly for $x>\sqrt t> 2$.
\end{Prop} 
\v2\n
\pf~  We can proceed as in the preceding proof except that in place of  (\ref{P_eq2})  we make decomposition
\[
x^{2\nu}\Big(G_{\la}(x) - (1/x)^{2\nu} G_{\la}(1) \Big) = G_{x^2\la} (1)-  G(1)+ [G(1)-G_{\la}(1) ]
\]
and  estimate the contributions of the three terms on the right hand side separately.  Let $\nu\neq 1$. It follows from the preceding   proof that the contributions of the last two terms
to $\eta(x,t)$ is bounded by a constant multiple of $x^{-2\nu}t^{-\nu-1}=t^{-2\nu-1}(\sqrt t\,/x)^{2\nu} $. As for the first term, on the one hand,  we recall that $G_z(1) = C_\nu(2z)^{\nu/2}K_\nu(\sqrt{2z})= C_\nu' z^{(2\nu -1)/4} e^{-\sqrt{2z}}(1+o(1))$ as $z\to\infty$ to see that 
$$x^{-2\nu}|f(-u\pm i0)G_{-x^2 u \pm i0} (1)|\leq  C u^{1\wedge \nu} x^{-2\nu}(x^2 u)^{(2\nu-1)/4} = C u^{1\wedge \nu + \nu} (x^2 u)^{-\nu/2-1/4}$$
for $u>1/x^2$ ($f$ is the same as before) and the integration over $u>1/x^2$ of $e^{-tu}$ times the last member yields the magnitude of the order $O\Big(t^{-(1\wedge \nu)+\nu+1}(\sqrt{t}/x)^{\nu+1/2}\Big)$. On the other hand we have
$x^{-2\nu}|f(-u\pm i0)G_{-x^2 u \pm i0} (1)|\leq C  u^{1\wedge \nu} x^{-2\nu}$  for $0<u<1/x^2$ and the corresponding integral  does not exceed the foregoing  magnitude.   By (\ref{P_eq0}) we find the asserted bound for $\nu\neq1$. The case $\nu=1$ is omitted, it being   similarly dealt with.
\v2\n
{\sc Remark about the case $\nu=0$ and proof of Corollary \ref{cor5}.} 

  Theorem \ref{thm1}  (as well as  any other results of \cite{Ubh})  does not give precise asymptotic form  for $x$ near $a$, while the first formula of Corollary  \ref{cor5} does. Here we indicate a manner by which  such an estimate can be obtained by following the method employed in the proof of Proposition \ref{P1}, and thereby prove the first formula  of Corollary \ref{cor5}.  By the way this entails the case  $x/\sqrt t \to 0$ of  the formula (\ref{R21}).

 For $\nu=0$ we must replace $G_z(x)/G_z(a)$ by $K_0(x\sqrt{2z})/K_0(a\sqrt{2z})$ in the inversion formula (\ref{P1.1}) that represents $p(x,t;a)$. Put $g(\la) =  -\lg (a\sqrt{\k \la})$ ($\k= 2 e^{-2\ga}$) as in \cite{Ubh};
 $g(\la)$ is the principal part of $K_0(a\sqrt{2\la})$ as $\la\downarrow 0$. The analogue of  the decomposition (\ref{P1.3}) should be
 $$\frac{K_0(x\sqrt{2\la})}{K_0(a\sqrt{2\la})} = \frac{K_0(x\sqrt{2\la})-K_0(a\sqrt{2\la})}{g(\la)} + R(\la;x),$$
 where
 $$R(\la;x) = \bigg[\frac1{K_0(a\sqrt{2\la})}-\frac1{g(\la)}\bigg]\bigg[ K_0(x\sqrt{2\la}) - K_0(a\sqrt{2\la}) \bigg] +1,$$
  With this $R$   define  $\eta$ by (\ref{P_eq0}) for which we have (\ref{P_eq5}).  Then we can proceed as in the proof of Theorem 1 of \cite{Ubh} with the same computation but taking  care of the effect of the subtraction of $K_0(a\sqrt{2\la})$  in the expressions above that results in the additional factor $(x-a)\wedge 1$ in front of the error term and  hence allows us to replace $1+\lg_+ x$ by $\lg (x/a)$ in the error estimate of   Theorem 1 of \cite{Ubh}  so that uniformly for  $x>a$, as $t\to\infty$
  \beqn\label{q_2}
  q^{(2)}(x,t;a)  = 2[\lg (x/a)] \frac{\k }{a^2} W\bigg(\frac{\k }{a^2}t\bigg) + O\bigg(\frac{\lg (x/a)}{t(\lg t)^2}\cdot\frac{x^2\wedge t}{t}\bigg),
  \eeqn
where $W(\la)= \int_0^\infty \frac{e^{-\la u}du}{[\lg u]^2 +\pi^2} = \frac1{\la(\lg \la)^2} -\frac{2\ga}{(\lg \la)^3} +\cdots.$
   
   Now let $a=1$ and substitute  $\frac12 \k =e^{ -2\ga}$ into the numerator of the leading term in the formula of Theorem \ref{thm1}, and   you find that for $1<x <\sqrt t$,
   $$q^{(2)}(x,t) =  \frac{2\lg x}{t[\lg(\k t)]^2}e^{-x^2/2t} +  \frac{2\ga (1- e^{-x^2/2t})}{t[\lg( \k t)]^2} +\frac{-4\ga\lg x}{t[\lg( \k t)]^3} +O\bigg(\frac{1}{t[\lg t]^3} \bigg).$$
   Using this estimate  for $\sqrt t/\lg t< x <\sqrt t$ and (\ref{q_2}) for  $1<x\leq \sqrt t/\lg t$ we obtain
  $$q^{(2)}(x,t) =  \frac{2\lg x}{t[\lg(\k t)]^2}e^{-x^2/2t}  +O\bigg(\frac{\lg  x}{t[\lg t]^3} \bigg)~~~~~~(a<x<\sqrt t)$$
  as $t\to\infty$, which is the same as  the first formula of Corollary \ref{cor5}.
  
  \section{Proof of Theorem \ref{thm3}}
For  the proof of  Theorem \ref{thm3} it suffices to verify  the formula of it in the case   $\sqrt t< x <Mt$ for each $M>1$ in view of Theorem \ref{thm3} and Remark given at the end of Section 2. It is convenient to  treat the cases $\nu=0$ and $\nu>0$ separately. In  both cases one may suppose  $x/t$ to tend to a constant $v\geq 0$ and  the subcases $v = 0$ and $v>0$ are also  separately treated since different arguments are employed for them, although the  framework  is the same.
In the case $v=0$ we shall  provide   estimates of the error terms   that are not given  in  Theorem \ref{thm3}.
\v2\n
{\bf 3.1. ~ The case $\nu=0$} 

In the course of proof of Theorem \ref{thm3} given below we shall derive
the following proposition, which   entails  the formula of  Theorem  \ref{thm3}   in the case $x/t \to 0$.
 \begin{Prop}\label{lem3.0}~   Let $\nu=0$.  It holds that uniformly for $ \sqrt t < x< t/2$,  as $t\to\infty$ 
 $$q^{(2)}(x,t, a)   =   \frac{\pi}{K_0(a x/t)}\, p^{(2)}_t(x) \bigg[1 + O\bigg( \frac{1}{\lg(t/ x)}\bigg) \bigg].$$
\end{Prop}

For the proof of Proposition  \ref{lem3.0} we shall use Theorem \ref{thm1}, which it is convenient to reduce  to the following slightly weaker  form. 
\begin{Lem}\label{lem3.1}~ Let $\nu=0$.
Uniformly for $x >a$, as $t\to\infty$
\beq
q^{(2)}(x,t; a)   &=& 2\pi p^{(2)}_t(x) \Bigg[  \frac{2\lg (x/a)\,}{[\,\lg ( t/a^2)]^2} + O\bigg( \frac{1}{ (\lg  t)^2} \bigg) \Bigg]~~~~\mbox{if}~~~~ x^2< 2t\lg (\lg t), \\
&=&   2\pi p^{(2)}_t(x) \Bigg[  \frac{2\lg (x/a)\,}{[\,\lg ( t/a^2)]^2} +  o\bigg( \frac{1}{ \lg  t} \bigg) \Bigg]  ~~~~~\mbox{if}~~~~2t\lg (\lg t)\leq x^2\leq 4t\lg (\lg t).
\eeq 
\end{Lem}
\v2\n

In what follows we let $a=1$. 
We need the following  lemma from \cite{Usaus} in which $\nu=0$. 
\begin{Lem}\label{lem3.11} ~ For any $\nu$, there  is a  constant $c=c_\nu>0$ such that for all $x>1$ and  $t>1$,
 \beqn\label{density_bound} 
q^{(d)}(x,t) \leq  c p^{(d)}_{t+1}(x). 
\eeqn 
\end{Lem}
\v2\n
\pf~ The proof follows from the parabolic Harnack inequality (cf. eg., {\cite{Ev}, \cite{W1})  as in the case $\nu=0$.
\v2

We use the fact that the Bessel process of order $\nu=0$ is the radial process of the standard two dimensional Brownian motion $B^{(2)}_t$. Let $\x$ denote a generic point of $\R^2$ and  $P^{BM(2)}_{\x}$ the probability of $B^{(2)}$ started at $\x$ . 
We can suppose that the initial point $B^{(2)}_0=\x$ is on the upper vertical axis so that $\x =(0,x)$. Write  $\xi_t$ and   $Y_t$  for the horizontal and vertical components of $ B^{(2)}_t$, respectively, and
  let $T_K$ be the first hitting time of the vertical level $K$ by $B^{(2)}_t$: $T_K=\inf\{t>0: Y_t=K\}$. Then the space-time distribution of  $(T_K, \xi_{T_K})$ is given by
\beqn\label{2_1}
\frac{P^{BM(2)}_{(0,x)}[T_K\in dt, \xi_{T_K}\in d\xi]}{dtd\xi} = \frac{x-K}{t}\,p^{(1)}_{t}( x-K)p^{(1)}_{t}( \xi)
\eeqn
(cf. \cite{IM}, page 25), 
which yields the representation
\beqn\label{f1}
q(x,t)= \int_0^t ds\int_{-\infty}^\infty \frac{x-K}{t-s}\, p^{(2)}_{t-s}\Big(\sqrt{\xi^2+( x-K)^2}\,\Big)q\Big(\sqrt{\xi^2+K^2}, s\Big)d\xi.
\eeqn
 $K$ may be any number between $1$ and $x=|\x|$ but we suppose $4\leq K <x/2$. With a fixed $K$ we  are to compute the repeated integral on the right-hand side by using the formula of  Lemma  \ref{lem3.1}. It is remarked that we make no use of  Lemma  \ref{lem3.1}  in the case when $x /t\to v >0$.

We break the rest of the proof  into four  parts. 
For the case $x/t\to 0$ certain  elaborate computations directly yield the desired formula of Proposition  \ref{lem3.0}  (Parts 1 and 2).  In the case  $x /t\to v \neq 0$,   on the other hand, we first show the existence of limit of the ratio $q(x,t)/p_t(x)$ (Part 3).  While  it is difficult to identify the limit along the same line as in the case $x/t\to 0$,   
with its existence at hand another way  determines the limit (Part 4).

Throughout  the proof we suppose that for some $M>0$,
$$\sqrt{t} < x< Mt.$$
The constant $K\geq 4 $ may  be  fixed arbitrarily prior to Part 3, in which  we  need to take $K$ large enough, so we do not assign   $K$  a specific value. We  write $p_{t}$, $q$ for  $p^{(2)}_{t}$, $q^{(2)}$  to be consistent to our convention that the super-script $\,^{(d)}$ is dropped, while we continue to write $p^{(\k)}_t$ if $\k\neq 2$.
We put,  for $0\leq c< \tau\leq t$,
$$I_{c,\tau}= I_{c,\tau}(x,t) := \int_c^{\tau} ds \int_{-\infty}^\infty   \frac{\, x-K\, }{t-s}p_{t-s}\Big(\sqrt{\xi^2+( x-K)^2}\,\Big) q\Big(\sqrt{\xi^2+K^2}\,, s\Big) d\xi,$$
the contribution to the integral in (\ref{f1}) from the interval $c< s<\tau$. 

\v2
{\it Part 1: Estimation of} $I_{c,t}$.   Here  $c$ is a constant  not less than   $4$. 
In the identity 
$
p\a^2+ q\b^2 =pq(\a-\b)^2+(p\a+ q\b)^2$, where  $p, q\in \R$ with $ p+q=1$ and 
$\a, \b$  may be vectors in any Euclidian space,
 take  $p=(t-s)/t$ and   divide  both sides of it  by 
$$T= pqt =\frac{s(t-s)}{t}$$
to obtain
\beqn
\label{id2}\frac{1}{s}\a^2 +\frac{1}{t-s}\b^2 =  \frac1{t} (\b-\a)^2 + \frac1{T}\bigg(\a +\frac{s}{t}(\b-\a)\bigg)^2.
\eeqn
Then substituting the two-dimensional vectors $\a=(\xi, K), \b =(\xi, K-x)$ leads to
\beqn\label{f2}
p_{t-s}\Big(\sqrt{\xi^2+( x-K)^2}\,\Big) p_s\Big(\sqrt{\xi^2+K^2}\,\Big)= \frac1{2\pi T}\,p_t(x) e^{-\xi^2/2T}\exp\Bigg[- \frac1{2T}\bigg(K-\frac{s}{t} x\bigg)^2\Bigg].
\eeqn

In the repeated integral of $I_{c,t}$ we  split the range of integration w.r.t. $\xi$ at $\xi =\pm \sqrt{4s \lg \lg s}$. We claim  that 
\begin{eqnarray}\label{f3}
I_{c, t} &=&p_t(x) \int_{c}^{ t} \frac{\, x-K\,}{(t-s)\sqrt T \,} \exp\Bigg\{- \frac1{2T}\bigg(K-\frac{s}{t}x\bigg)^2\Bigg\} \times \nonumber\\
&&~~~~~~\times \Bigg[ \int_{\sqrt{\xi^2 +K^2} <\sqrt{4 s\lg\lg s}} \frac{\,\lg (\xi^2+K^2)\,}{(\lg s)^2\sqrt T} e^{-\xi^2/2T} d\xi  \nonumber\\
&&~~~~~~~~~~~~~~~~~~~~~~~~~~ + R(s,t)+ O\bigg( \frac{1}{(\lg s)^2} \bigg)\Bigg]ds,
 \end{eqnarray}
where $R(s,t)$, the term that comes from the remainder term in Lemma  \ref{lem3.1}, is   $o(1/\lg s)$. 
For the part $|\xi|\geq \sqrt{4 s\lg s}$ we have only to substitute  the expression  of $q$ given in Lemma  \ref{lem3.1} and apply (\ref{f2}) (note that the bound in Lemma \ref{lem3.1} actually holds uniformly for $t\geq 4$ simply because $q(x,t)$ is bounded there).
For the integral on $|\xi|\geq \sqrt{4 s\lg s}$ we need to take  $p_{s+1}$ in place of $p_s$ in (\ref{f2}) so that the corresponding contribution to  $I_{c, t}$ becomes
\begin{eqnarray}\label{f300}
&&p_{t+1}(x) \int_{c}^{ t} \frac{\, x-K\,}{(t-s)\sqrt {T' }} \exp\Bigg\{- \frac1{2T'}\bigg(K-\frac{s+1}{t+1}x\bigg)^2\Bigg\}\times \nonumber\\
&&~~~~~~~~~~~~~~~~~~\times  \int_{\sqrt{\xi^2 +K^2} \geq\sqrt{4 s\lg\lg s}} \frac{\,e^{-\xi^2/2T'} \,}{\sqrt {T'}}  d\xi ds,
 \end{eqnarray}
where $T' = (s+1)(t-s)/(t+1)$; 
 since the inner integral  is $O(1/(\lg s)^2)$ uniformly in $t$,  a simple change of variable  gives the error term in (\ref{f3}). 
 Thus we have verified  the claim.
  
  Scaling the variable $\xi$ by $\sqrt s$  and dominating     $e^{-\xi^2/2T}$ by  $e^{-\xi^2/2s}$ if necessary  we see that the quantity in the big square brackets is evaluated to be   
\beqn\label{f4}
\frac{\sqrt{2\pi} \,}{ \lg s}  +O\bigg(\frac{1  }{(\lg s)^2\sqrt{(t-s)/t}\,}\bigg) +R(s,t).
\eeqn
 We must  compute the integral
\beqn\label{f41}
 J:=  \int_{c}^{ t} \frac{\, x-K\,}{(t-s)\sqrt T \,} \exp\Bigg\{- \frac1{2T}\bigg(K-\frac{s}{t}x\bigg)^2\Bigg\} \frac{\sqrt{2\pi} \,}{\lg s}ds.
\eeqn

Now we suppose $x/t <1/2$ so that $\lg t/x >\lg 2$ and   take  $c=4$.  By a simple change of  the variables of integration, e.g.,  according to $u= (x/t) \sqrt s$, which transforms  $s/t $  to $u^2t/x^2$  one can easily find that  $J \sim  p_t(x)\pi/\lg(t/x)$ as $x/t \to 0$,  $t/x^2 \to 0$  (which is enough for Theorem \ref{thm3} restricted to the case $x/t \to 0$). But this way does not give the error estimate  asserted in Proposition \ref{lem3.0}. To improve the evaluation of the integral
we transform  the  variable of integration by
\beqn\label{c_v}
\rho =\frac{s}{t-s}~~~~~~~~\bigg(\mbox{i.e., }~~ s= \frac{t\rho}{1+\rho}\bigg),
\eeqn
entailing the relation $ds= (t-s)^2 d\rho/t = (t-s) \sqrt {T/t\rho} \, d\rho$. 

Noting the inequalities 
$$0<\frac1{T} -\frac{1}{t\rho} \le   \frac2{t-s}~~~~\mbox{ and}~~~~  \frac1{T}\cdot \frac{sx}{t}=\frac{x}{t-s}<\frac{2x}{t}\vee \bigg(\frac2{x}\cdot \frac{x^2\rho}{t}\bigg)$$
we may write the exponent appearing in the integral  of (\ref{f41})  in the form
\beq\label{f5}
 - \frac1{2T}\bigg(K-\frac{s}{t}x\bigg)^2 = -\frac{K^2}{2t\rho} - (1 -\de)\frac{ x^2 \rho}{2t}  +O\bigg(\frac{x}{t}\bigg),
\eeq  
where $\de=\de(t,x,s, K)$ is a function of $t,x, s, K$ that satisfies $0<\de <4K/x $ (provided that $K<x$).
Further transform  the variable  $\rho$  to  $u =x\sqrt {\rho/t}$. 
Then
\beq\label{c_v2}
\frac{xds}{(t-s)\sqrt T\,} =\frac{x d\rho}{\sqrt{t\rho}} =2du,
\eeq
and we obtain    
\beqn\label{upperlimit}
 J = \int_{(x/t)\sqrt{c/(1-c/t)}}^{\infty} \frac{2\sqrt{2\pi}}{\, \lg \Big[(t/x)^2 m(u)
\Big]\,}\exp\bigg( -\frac{K^2x^2}{2t^2 u^2} - (1-\de)\frac{u^2}{2}\bigg)du \,\Bigg(1+O\bigg(\frac{x}{t}\bigg)\Bigg), 
\eeqn
where $m(u) =u^2/(1+u^2t/x^2)$. We apply 
the inequality  $ |1-1/(1+r)| \leq   |r|+ r^2/(1+r)$ ($r>-1$) with   $r=[\lg m(u)]/\lg(t^2/x^2)$ for which $(1+r)^{-1}\leq (\lg c)^{-1}\lg(t^2/x^2)$ to see that
$$ \frac{2\sqrt{2\pi}}{ \lg \,[(t/x)^2 m(u)]} =\frac{\sqrt{2\pi}}{\, \lg (t/x)} + O\bigg(\frac{1+ [\lg m(u)]^2}{[\lg (t/x)]^2}\bigg)$$
uniformly valid if $u$ is confined to  the range of integration.   
Owing to  the identity 
\beqn\label{exp_int}
\int_0^\infty e^{- b/2u^2 - \la u^2 /2 }du =\sqrt{\frac{\pi}{2\la \,} } \,e^{-\sqrt{2 b\la}}~~~~~~~~~~~ (\la>0, b\geq 0)
\eeqn
and the bound  $\int_0^\infty |\lg m(u)|^2e^{- u^2/2}du\leq C$  the formula (\ref{upperlimit}) reduces  to
\beqn\label{conc}
J =[\pi /\lg (t/x)] [1+ O(1/\lg (t/x) )] .
\eeqn
Taking the computation carried out right above  into account one also observes that
the contribution of  the error   term  in (\ref{f4})  is  $O\Big(p_t(x)/[\lg (t/x)]^2\Big)$ and  concludes that  uniformly for $\sqrt t<x<t/2$, as $t\to\infty$
\beqn\label{R0}
I_{c,t} = p_t(x)\Bigg[ \frac{\pi}{\, \lg ( t/x)\,} \Bigg(1+   O\bigg(\frac1{\lg (t/x)}\bigg) + \tilde R\Bigg)\Bigg],
\eeqn
 where  $\tilde R$, the term corresponding to $R(s,t)$, is $o(1)$.  
\v2
{\it Part 2: Estimation of} $I_{0,c}$. ~  
The integral $I_{0,c}(t,x)$ is dominated by
\beq
&&\int_{-\infty}^\infty P\!\!_{\sqrt{\xi^2+K^2}}\, [\sigma_1 < c] \sup_{0< s<c}  \frac{x-K}{t-s}\,p_{t-s}\Big(\sqrt{\xi^2+( x-K)^2}\,\Big)d\xi \\
&& ~~~~~~\leq C\sqrt{c}  \,e^{-K^2/2c} \frac{x}{t K} p_t( x-K)  
\leq \Big[C' e^{-K^2/2c}e^{Kx/t}\Big] p_t(x) \frac{x}{t},
\eeq
where for the first inequality  we have used the bound $P_y[\sigma_1<c]\leq C\sqrt{c} e^{-(y-1)^2/2c}/(y-1)$ for $y>1$, a bound obtained from the one-dimensional result (cf. Lemma 3.2 of \cite{Usaus}). 
Combined with (\ref{R0})  this shows that $q(x,t) /p_t(x)= O(1/\lg t)$ at least for $\sqrt{2 t\lg\lg t}<x \leq \sqrt{4 t\lg\lg t}$. Using this bound instead of the second one of Lemma \ref{lem3.1} we obtain
$$R(s,t) = O(1/(\lg s)^2)$$
so that $\tilde R = O(1/\lg(t/x))$.
The proof of Proposition \ref{lem3.0} is complete.

\v2
\v2
{\it Part 3: Proof of convergence. }~ Here we suppose  $x/ t \to v >0$ and  prove that there exists $\lim q(x,t)/p_t(x)$,  the limit value depends only on $v$ and  the  convergence is locally uniform in $v$.  Here we use Lemma \ref{lem3.11} but does not  Lemma \ref{lem3.1}. With the help of (\ref{f2}) it gives   
\beqn\label{f4_0}
 I_{K^2, \,t/2} (x,  t)\leq Cp_t(x) \int_{K^2}^{ t/2} \frac{\, x\,}{t} \exp\Bigg[- \frac1{2(1-s'/t')}\bigg(\frac{K}{\sqrt {s'}}-   \frac{\, x\,}{t'}{\sqrt {s'} } \bigg)^2\Bigg]  \frac{\,ds\,}{\sqrt s },
\eeqn
where  $s'= s+1$, $t'=t+1$ and (\ref{f2}) is applied with  $p_{s+1}\Big(\sqrt{\xi^2+K^2}\,\Big)$ in place of  $p_{s}\Big(\sqrt{\xi^2+K^2}\,\Big)$ (see (\ref{f300})).    One observes that  the integral above   is at most $O( e^{- v K/4})$ (use e.g. (\ref{exp_int})).  Also a quite crude estimation shows $ I_{t/2, t} (x,  t) \leq C p_t(x)e^{-(2K-x)^2/8t}$. 
Combined with 
   the result of  Part  2 these show  that  for any $\e>0$ one can choose  $K$ large so that
\beqn\label{conv}
\limsup_{ t\to\infty, x /t \to v}\Bigg| \frac{q(x,t) -I_{c, K^2}(x,t) }{p_t(x)}\Bigg| <\e.
\eeqn

Define $h_K(\xi, s)$ by 
$$h_K(\xi,s) =q\Big(\sqrt{\xi^2+K^2}\,,s\Big)\Big/p_s\Big(\sqrt{\xi^2+K^2}\,\Big).$$ 
By Lemma \ref{lem3.11} 
$$h_K(\xi,s) \leq C \exp\bigg[\frac{\xi^2+K^2}{2s(s+1)}\bigg]$$
 and keeping this bound in mind we see that
\beq
&&\frac{I_{c, K^2}(x,t) }{p_t(x)}   \\&&=  \frac1{p_t(x)}\,  \int_c^{K^2} ds \int_{-\infty}^\infty   \frac{\, x-K\,}{t-s}\,p_{t-s}\Big(\sqrt{\xi^2+(x-K)^2}\,\Big) p_s\Big(\sqrt{\xi^2+K^2}\,\Big)h_K(\xi, s)d\xi \\
&&=  \int_c^{K^2}\frac {ds}{2\pi T}   \int_{-\infty}^\infty  \frac{\, x-K\,}{t-s}\, \exp\Bigg[- \frac{1}{2T}\bigg( \Big(K- \frac{x}{t}s \Big)^2 +\xi^2\bigg)\Bigg] h_K(\xi, s) d\xi \\
&& \longrightarrow  ~ v  \int_c^{K^2}\frac{ ds}{2\pi s} \int_{-\infty}^\infty  \exp\Bigg[- \frac{\,(K-v s)^2\,}{2s} -\frac{\xi^2}{\,2s\,} \Bigg] h_K(\xi, s) d\xi
\eeq
as $x /t \to v$.
This together with (\ref{conv}) shows that $q(x,t)/p_t(x)$ is convergent  and the limit value does not depend on the manner of  $x/t$ approaching to $v$. The required  uniformity   of the convergence is easily ascertained from the arguments made above.

\v2
{\it Part 4: Identification of the limit.}~
 The proof rests on the identity
\beqn\label{Fun_eq}
p_t(x)= \int_0^tq(x, t-s)p_s(1)ds.
\eeqn
which follows from the Markov property of the Bessel process and also from the identity (\ref{LT}).
 We may suppose that $x = tv $, $v\neq 0$.  By Part 3
\beqn\label{f6}
q(t v, t-s) = \la p_{t-s}(t v)(1+o(1)) ~~~ \mbox{as} ~~ s/t \to 0, t\to\infty
\eeqn
for some  constant $\la =\la(v) \geq 0$.  
Since
$$\frac{p_{t-s}(t v)p_s(1)}{p_t(t v)} =  \frac{1}{2\pi s(1-s/t)}\exp\bigg(- \frac{v^2 s}{2(1-s/t)} - \frac{1}{2s}\bigg) $$
Substitution of (\ref{f6})  into (\ref{Fun_eq}) yields
\begin{eqnarray}\label{f7}
\frac{1}{\la}&=& \lim \frac1{p_t(t v)} \int_0^t p_{t-s}(tv))p_s(1) ds \\
&=& \frac1{2\pi}\int_0^\infty \exp\bigg(- \frac{v^2 s}{2} - \frac{1}{2s}\bigg) \frac{ds}{s} \nonumber\\
&=& K_0(v)/\pi    \nonumber
\end{eqnarray}
(see \cite{E}, Eq(29) in page 146 for the last equality).  Hence $\la = \pi/K_0(v)$ as desired.  This completes the proof of Theorem \ref{thm3} in the case when $x/t$ is bounded. \qed
\vskip2mm 
  {\sc Remark} 4.  In the case $x /t  \to v>0 $, it seems hard  to compute the value  $\lim q(x,t)/p_t(x)$ along the same line as in Part 1  since our knowledge of   the   behavior of $q((\xi, K), s)$ is poor for small values of $s$ that significantly  contributes  to the integral of (\ref{f1}). 
  This point would be well  understood from the argument made in Part 3 above. One notes that from  Part 2 we know only that $I_{0,c}(x,t)$ becomes small if $K/c$ is large enough, while  the error term  $O(1/(\lg T)^2)$ in  the estimate  (\ref{f3}) depends on  $c$ .
  \v2\n
{\bf 3.2. ~ The case $\nu>0.$}
  \begin{Prop}\label{thm3.2}~   Let $\nu>0$.  It holds that   uniformly for $ \sqrt t < x< t/2$,  as $t\to\infty$ 
 \beq
 \frac{q^{(d)}(x,t, a)}{p^{(d)}_t(x) }   -   \frac{1}{G(a)}\, &=&  O\Bigg(\bigg(\frac{x}{t}\bigg)^{d-2}\sqrt{\lg \frac{t}{ x}}\,\Bigg) ~~~~~(0<\nu<1/2),\\
 &=& O(x/t)~~~~~~~~~~~~~~~~~~~~~~~~(\nu>1/2).
 \eeq
\end{Prop}

 We use the fact that the square of Bessel process of dimension  $d>2$ is the sum of those of  two independent Bessel processes of dimensions 1 and $d' = d-1$  (\cite{RY}). Let  $(Y_t)$  be the one-dimensional standard Brownian motion started at $x$ and $ (\xi_t)$   the Bessel process of the dimension $d'$ started at $0$ and independent of $(Y_t)$ . Then
the law of the process $(X_t^2)$ is the same as the law of $(\xi_t^2 + Y_t^2)$. 
  Let $T_K=\inf\{t>0: Y_t=K\}$. Then in place of (\ref{2_1}) we have
\beq
\frac{P[T_K\in dt, \xi_{T_K}\in d\xi]}{dtd\xi} = \frac{x-K}{t}\,p^{(1)}_{t}( x-K)p^{(d')}_{t}( \xi)  c_{d'}\xi^{d-2}d\xi
\eeq
so that
\beqn\label{f1_0}
q(x,t)= \int_0^t ds\int_{-\infty}^\infty \frac{x-K}{t-s}\, p^{(1)}_{t-s }( x-K)p^{(d')}_{t-s}( \xi) q\Big(\sqrt{\xi^2+K^2}, s\Big) c_{d'}\xi^{d-2}d\xi.
\eeqn

The proof of Proposition  \ref{thm3.2} given below  is  analogous to the one given for $\nu=0$ and  we  proceed parallel to the lines of the  preceding proof.

\v2
{\it Part 1: Estimation of} $I_{c, t}$.    We write $I_{b,c}(x,t)$ as before for  the integral in (\ref{f1_0}) restricted on
the interval $[b,c]$. 
 The product $ p^{(1)}_{t-s }( x-K)p^{(d')}_{t-s}( \xi)$  appearing in the integrand may be written as
\beq
p_{t-s}\Big(\sqrt{\xi^2+( x-K)^2}\,\Big) =  
 p_{t-s}(x-K)e^{-\xi^2/2(t-s)}.
\eeq
 which we further rewrite in the form
\beqn\label{two_id}
\bigg(\frac{t}{t-s}\bigg)^{d/2} p_t(x-K) \exp\bigg\{- \frac{(x-K)^2s}{2t(t-s)} \bigg\} e^{-\xi^2/2(t-s)}.
\eeqn

We split the range of $\xi$-integration at $\xi =\pm \sqrt{4s \lg s}$  in the repeated  integral of  $I_{c, t}$.  The integral on $|\xi|\geq \sqrt{4 s\lg s}$ is disposed of  by employing Lemma \ref{lem3.11}  as before (see (\ref{f300})). 
As for  the integral on the other part
 we first  evaluate the contribution of the term $(\xi^2+K^2)^{-\nu}p_{s}(1)/G(1)$, which, on using  (\ref{two_id}),    is dominated  by a constant multiple of
$$ R_1:= p_t(x)\int_{c}^{ t} ds   \frac{\, xt^{d/2} p_s(1)\, }{(t-s)^{d/2+1}}\exp\bigg\{- \frac{(x-K)^2s}{2t(t-s)}  \bigg\} \int_{|\xi|< \sqrt{4s \lg s}} e^{-\xi^2/2(t-s)}  d\xi. $$
It is convenient to split the outer integral at $s=t/2$ and let $R_{1}^>$ and $R_{1}^<$ be the parts corresponding to $s>t/2$ and $s\leq t/2$, respectively. By performing the $\xi$-integration  and changing the variable by $u=t-s$ we deduce 
$$R_{1}^> \leq  C p_t(x) \int_0^{t/2} \frac{x}{u^{(d+1)/2}} e^{- x^2/4u}du\leq C'' p_t(x) x^{2-d}.$$
For the evaluation of $R_1^<$
we replace the integrand by unity in the inner integral  and  have the bound
\[
R_1^<\leq  \frac{Cp_t(x)x}{t} \int_{c}^{ t/2}  \frac{ \sqrt{\lg s}}{ s^{(d-1)/2}} \exp\bigg\{- \frac{x^2s}{2t^2}\bigg\} ds. 
\]
Since the integral on the right-hand side is evaluated to be $O\Big(x/t)^{d-3}\sqrt{\lg t/x}\Big) $ or $O(1)$ according as $\nu< 1/2$ or $\nu>1/2$, by taking  account of the estimate for $R_1^>$ obtained above we deduce that
\beqn\label{R_1}
R_1 \leq C' p_t(x) \bigg(\frac{x}{t}\bigg)^{d-2}\sqrt{\lg \frac{t}{x}} ~~~\mbox{if ~$\nu<\frac12$~~~and~~~}
R_1 \leq C' p_t(x) \frac{x}{t} ~~~\mbox{if ~$\nu>\frac12$}.
\eeqn

Let $0<\nu <1/2$. Then, in a similar way to the above, we evaluate  the contribution of the error term in  (\ref{R1b}), denoted by $R_2$ and make decomposition $R_2 =R_2^> +R_2^<$. For $R_2^>$ we note that $\int_\R e^{-\xi^2/2(t-s)}\xi^{d-2}d \xi = O((t-s)^{(d-1)/2})$ and deduce that 
$$
|R_2^>|\leq Cp_t(x) xt^{-\nu} \int_0^{t/2} u^{-3/2}e^{-x^2/4u}du\leq  C'p_t(x) t^{-\nu};$$
also
\begin{eqnarray}\label{R_222}
|R_2^<| &\leq&  \frac{Cp_t(x)x}{t} \int_{c}^{ t/2}  
   \frac{e^{-(x^2/2t^2)s}ds}{s^{d/2 +\nu}} \int_0^{\sqrt{4s\lg s}}  \xi^{d-2}d\xi  \nonumber\\
   &=&\frac{Cp_t(x)x}{t} \int_{c}^{ t}  
   \frac{(\lg s)^{(d-1)/2}}{s^{(d-1)/2}} e^{-(x^2/2t^2)s}ds
 \leq C' p_t(x) \bigg(\frac{x}{t}\bigg)^{d-2}\sqrt{\lg \frac{t}{x}},
\end{eqnarray}
so that
\beqn\label{R_2}
|R_2|  \leq C''' p_t(x)\bigg(\frac{x}{t}\bigg)^{d-2}\sqrt{\lg \frac{t}{x}}.
\eeqn

Putting $T= s(t-s)/t$ we have in place of (\ref{f2})
\beq\label{f2_0}
p_{t-s}\Big(\sqrt{\xi^2+( x-K)^2}\,\Big) p_s\Big(\sqrt{\xi^2+K^2}\,\Big)= p_t(x) p_T^{(d')}(\xi)\frac1{\sqrt {2\pi T}}\exp\Bigg[- \frac1{2T}\bigg(K-\frac{s}{t} x\bigg)^2\Bigg].
\eeq
 Applying  this together with  (\ref{R_1})  and  (\ref{R_2})  and making use of Lemma \ref{lem3.11}  in the same manner as before we find that 
\begin{eqnarray*}\label{f3_0}
I_{c, t} &=&\frac{p_t(x)}{G(1)} \int_{c}^{ t} \frac{\, x-K\,}{(t-s)\sqrt{2\pi T} \,} \exp\Bigg\{- \frac1{2T}\bigg(K-\frac{s}{t}x\bigg)^2\Bigg\} \times \nonumber\\
&&~~~~~~\times \Bigg[ \int_{\sqrt{\xi^2 +K^2} <\sqrt{4 s\lg s}} p_T^{(d')}(\xi) c_{d'}\xi^{d-2}d\xi  + O\bigg( \frac{1}{s} \bigg)\Bigg]ds \\
&&~~~~~~~~~+ O\Bigg( p_t(x) \bigg(\frac{x}{t}\bigg)^{d-2}\sqrt{\lg \frac{t}{x}}  \,\Bigg).
 \end{eqnarray*}
The quantity in the big square brackets may be evaluated to be   $1+O\Big(1/s[(t-s)/t]^{(d-1)/2}\Big)$.  In order to evaluate the whole integral  we employ the transformation (\ref{c_v})  and  follow the succeeding  arguments up to  (\ref{conc}). We then conclude that 
$$
I_{c, t}  =   \frac{p_t(x)}{G(1)} \Bigg[ 1+  O\Bigg( \bigg(\frac{x}{t}\bigg)^{d-2}\sqrt{\lg \frac{t}{x}}  \,\Bigg) \Bigg].   
$$

Let  $\nu>1/2$. Then,  the  integral of the third member  in  (\ref{R_222})  becomes bounded, so that we have  $|R_2^<| \leq C'p_t(x) x/t$ in place of the bound given therein.    The other computations may be carried out in a similar way and we obtain $I_{c,t}=[p_t(x)/G(1)] (1+O(x/t))$.

\v2
{\it Part 2: Estimation of} $I_{0,c}$. ~
The same computation as before gives  the same bound of $I_{0,c}$ (but here $p_t= p_t^{(d)}$), which is sufficient for the present purpose.   Thus Proposition \ref{thm3.2} has been proved.

\v2

\v2
{\it Part 3: Proof of convergence.}~ 
The proof is quite similar to the one given for $\nu=0$.  The bound  (\ref{f4_0}) and  hence the relation  (\ref{conv}) holds true without any alteration except that  here $q$ and $p_t$ are defined with   $d>2.$
Define $h_K(\xi, s)$ as before. Then 
\[
\frac{I_{c, K^2}(x,t) }{p_t(x)} ~\longrightarrow  ~ v  \int_a^{K^2}\frac{ ds}{(2\pi s)^{d/2}} \int_{-\infty}^\infty  \exp\Bigg[- \frac{\,(K-v s)^2\,}{2s} -\frac{\xi^2}{\,2s\,} \Bigg] h_K(\xi, s) c_{d'}\xi^{d-2}d\xi
\]
as $x /t \to v$
and as before we conclude the desired convergence.
\v2
{\it Part 4}.~
  Let  $\la_d(v)$ be the limit of $q(x,t)/ p_t(x)$ as $x/t \to v > 0$.  The functional equation (\ref{Fun_eq})  holds  true for all $\nu>0$  in view of the corresponding identity for the Laplace transforms. In place of (\ref{f7}) we then have that  if $x/t\to v$, then
\begin{eqnarray*}\label{f8}
\frac{1}{\la_d(v)}= \lim  \int_0^t \frac{p_{t-s}(tv)p_s(1) } {p_t(tv)} ds 
&=& \frac1{(2\pi)^{d/2}}\int_0^\infty \exp\bigg(- \frac{v^2 s}{2} - \frac{1}{2s}\bigg) \frac{ds}{s^{d/2}} \nonumber\\
&=&2v^{d/2-1} K_{d/2-1}(v)/(2\pi)^{d/2},    \nonumber
\end{eqnarray*}
so that $\la_d(v)= \La_d(v)$. Thus we  conclude the formula of Theorem \ref{thm3}. \qed

\section{Asymptotics of the distribution of $\sigma_a$}
We  derive estimates of  the distribution   $P_x[\sigma_a<t]$ from those of the density.  Here we compute only the 
principal parts of $P_x[\sigma_a<t]$ or $P_x[t<\sigma_a <\infty]$ (according as $t<x^2$ or $t\geq x^2$). With a little more labore one can obtain the error term by employing Propositions \ref{P1}, \ref{P2}, \ref{lem3.0}, \ref{thm3.2} or  Lemma \ref{thm5}. 
Let $\ga(y,\nu)$ be the incomplete gamma function and put
$$\ga_\nu(y) =  \frac{\ga(y,\nu)}{\Ga(\nu)} = \frac{1}{\Ga(\nu)}\int_0^{y}  e^{-u}u^{\nu-1}du.$$ 
and
$$\Ga_\nu(y) =  1- \ga_\nu(y) = \frac{1}{\Ga(\nu)}\int_y^{\infty}  e^{-u}u^{\nu-1}du.$$

\begin{Thm}\label{thm6}~ Let $\nu > 0$. Uniformly for  $x > a$, as $t\to \infty$, 
\beqn\label{X1}
\frac{P_x[t< \sigma_a<\infty]}{P_x[\sigma_a<\infty] }=\bigg[1-\bigg(\frac{a}{x}\bigg)^{2\nu}\, \bigg]\ga_\nu\bigg(\frac{x^2}{2t}\bigg)(1+o(1))
\eeqn
and
\beqn\label{X2}
\frac{P_x[\sigma_a<t]}{P_x[\sigma_a <\infty]} = \frac1{\La_\nu(0)}\La_{\nu}\bigg(\frac{ax}{t}\bigg)\Ga_\nu\bigg(\frac{x^2}{2t}\bigg)(1+o(1)).
\eeqn
\end{Thm}
\v2

If $x^2/t \to \infty$, the first formula (\ref{X1})     is poor (at least in comparison  with the second one)  since then $\ga_\nu(x^2/2t)$ tends to 1 and  it says simply that  the conditional probability of escaping from $a$ after $t$  tends to 1 and nothing more: such a result may be verified more directly (e.g.,    the crude bound  given in Lemma \ref{lem3.11} may be used to derive  a fairly nice estimate of the speed of convergence).   Similarly,  in the case $x^2/t \to 0$, (\ref{X2}) asserts that the conditional probability of  arriving $a$ before $t$  tends to 1, which also readily follows from
 (\ref{EQ11}) as well as from (\ref{X1}).

Taking limit along $x^2/2t = 1/y$, either of  (\ref{X1}) or (\ref{X2}) shows that the scaled variable $2\sigma_a/x^2$  conditioned on the event $\sigma_a<\infty$ converges in  law 
to  a variable whose  distribution function is   $1- \ga_\nu(1/y)$. This, however, follows immediately from (\ref{T(2)})  by knowing the formulae $K_\nu(t)\sim 2^{\nu-1}\Ga(\nu)t^{-\nu}$ ($t\downarrow 0$) (\cite{L}: page 136) and  $-\int_0^\infty e^{-\la y}d\ga_\nu(1/y) =  2K_\nu(2\sqrt \la\,)\la^{\nu/2}/\Ga(\nu)$ (\cite{E}:(29) on page 146).

The proofs of two formulae of Theorem \ref{thm6} are similar. Since  (\ref{X1})  is easier we prove  only  (\ref{X2}). By what is remarked above we have only to prove  it  for $x> \sqrt{t/\lg t}$. With this restriction we can include the case $\nu=0$. We remind the readers  that
$$P_x[\sigma_a <\infty]  =\bigg(\frac{a}{x}\bigg)^{2\nu}.$$

\begin{Thm}\label{thm7}~ Let $\nu\geq 0$. Uniformly for  $x > \sqrt{t/\lg t}$, as $t\to \infty$, 
$$P_x[\sigma_a<t] = \La_{\nu}\bigg(\frac{ax}{t}\bigg)\bigg(\frac{a}{x}\bigg)^{2\nu}\frac{2^\nu }{(2\pi )^{\nu+1}}\int_{x^2/2t}^\infty e^{-y}y^{\nu-1}dy(1+o(1)).$$
\end{Thm}
\v2\pf~ Employing Lemma \ref{thm5} (if necessary) as well as Theorem \ref{thm3} and recalling $\La_\nu(y) = Cy^{-\nu+1/2}e^y(1+o(1))$ for $y>1$  one observes first that
$P_x[\sigma_a <\sqrt t]$ is negligible and then that
$$P_x[\sigma_a<t]  =a^{2\nu}\int_0^t \La_{\nu}\bigg(\frac{ax}{s}\bigg) p_s^{(d)}(x)ds(1+o(1)).$$
By a simple change of the variable the right-hand side  is transformed into
$$\bigg(\frac{a}{x}\bigg)^{2\nu}\frac1{(2\pi)^{d/2}}\int_{x^2/t}^\infty \La_{\nu}\bigg(\frac{ay}{x}\bigg) e^{-y/2} y^{d/2-2}dy (1+o(1)).$$
If $\nu>0$, the proof is easy from this expression and  the following argument is made to include the case $\nu=0$.  If $ \sqrt {t/\lg t} < x<\sqrt {t\lg t}$ (for instance), then $x^2/t>1/\lg t$ and $(\lg t)^2x/t \to 0$, and  hence  the range of  integration  may be restricted to the interval $[x^2/t,  (\lg t)^2x^2/t]$  in which $\La_\nu(ay/x)= \La_\nu(ax^2/t)(1+o(1))$ so that one may replace $\La_\nu(ay/x)$  by $\La_\nu(ax^2/t)$, yielding the desired formula after a simple change of the variable of integration.  The other case may be dealt with in a similar manner.  If $\sqrt{t\lg t} \leq  x <t$, then $x^2/t$ goes to infinity so that 
the upper limit of the integral may be replaced by $(1+\de)x^2/t$ for any $\de>0$  and  the required relation is reduced to $\La((1+\de)ax/t)/\La_\nu(ax/t) \to 1$ as $\de\downarrow 0$ uniformly in this region, which is plainly true. As for the case $x \geq t$ one has only to replace $\de$ by  $K/x$ with large $K$ and argue  analogously. The proof of Theorem \ref{thm7} is complete. \qed

In the case when $\nu=0$ and  $x<\sqrt{2 t\lg \lg t}$ a  precise asymptotic form    is obtained in \cite{Ubh}. Combined with  it  Theorem \ref{thm7} shows 
\begin{Cor}\label{cor6} ~ Let $\nu=0$.   Uniformly for  $x > 1$, as $t\to \infty$, 
\begin{eqnarray}\label{Y1}
P_x[\sigma_1<t] &=& 1-\frac{2\lg x}{\lg t}\bigg[1 - \frac{\lg (2e^{-\ga})}{\lg t}+ O\bigg(\frac{\frac1{\lg t}\vee \frac1{t} x^2}{\lg t} \bigg)\bigg] \quad      \mbox{for} \quad x<\sqrt t\\
&=& \frac1{2K_{0}(x/t)}\int_{x^2/2t}^\infty \frac{e^{-y}}{y} dy\,  (1+o(1))  \quad      \mbox{for} \quad x>\sqrt {t/\lg t}.
\label{Y2}
\end{eqnarray}
\end{Cor}
\v2

From Corollary \ref{cor6}  it follows  that if $\nu=0$ and  $x= \mu t^\a$ with  $\mu>0$ fixed, then as $t\to\infty$,
$$P_x[ \sigma_1 <t ]  ~~\longrightarrow ~~ (1- 2\a) \quad \mbox{if}\quad  0\leq  \a<1/2,$$
 $$P_x[ \sigma_1 <t ]  \,\sim \, \int_{\frac12 \mu^2t^{2\a-1}}^\infty  \frac{e^{-y}}{2y} dy \times
 \left\{ \begin{array} {ll}  
  {\displaystyle \frac1{(1-\a)\lg t}  } \quad&\mbox{if} \quad  1/2\leq \a <1\\ [2mm]
 {\displaystyle \frac1{K_0(\mu) } }\quad&\mbox{if} \quad  \a =1\\ [2mm]
 {\displaystyle \Big(\pi^{-1}2\mu t^{\a-1} \Big)^{1/2}e^{\mu t^{\a-1} } } \quad&\mbox{if} \quad   \a >1. 
  \end{array} \right.
 $$  
 \v2\n
In view of the identity $P_x[\sigma_a<t] = P_{x/a}[\sigma_1 <t/a^2]$,   the case $\a=1/2$ of this formula  is the same as (1.6) of  Spitzer \cite{S}, where it  is  used to derive his well-known test for a parabolic thinness  at infinity of space-time boundaries. By the same token an equivalent form in the case $0<\a<1/2$  is  $\lim_{a\downarrow 0} P_1[\lg \sigma_a\leq \ga \lg a^{-1} ] = \ga(2+\ga)^{-1}$ ($\ga= \a^{-1}-2$), which   is found  in \cite{IM}, Problem 4.6.4 but in terms of Laplace transform.

\vskip4mm
{\bf Acknowledgments.}~  I wish  to thank an  anonymous  referee for informing  me  of  the paper  \cite{BMR}  and   providing several valuable  comments  on    comparison of the results of  the original version of  the present paper with those of \cite{BMR}.

\end{document}